\documentclass[12pt]{amsart}

\usepackage{amsthm, amssymb}
\usepackage{amsfonts, amscd}
\usepackage{epsfig,multicol}
\theoremstyle{plain} \numberwithin{equation}{section}
\newtheorem{thm}{Theorem}[section]
\newtheorem{cor}[thm]{Corollary}
\newtheorem{lem}{Lemma}[section]
\newtheorem{prop}[thm]{Proposition}

\def \Hom {{\rm Hom}}
\def \Span {{\rm Span}}
\theoremstyle{remark}
\newtheorem{rem}{Remark}
\theoremstyle{example}
\newtheorem{exmp}{Example}

\begin{document}
\title[Graphs of 2-torus actions]{Graphs of 2-torus actions}
\author[Zhi L\"u]{Zhi L\"u}
\thanks{Supported by grants from NSFC (No. 10371020 and No. 10671034)}

\address{Institute of Mathematics, School of Mathematical Science, Fudan University,\newline Shanghai,
200433, P.R. China.} \email{zlu@fudan.edu.cn} \maketitle
\date{}
\begin{abstract} It has been known that an effective smooth
$({\Bbb Z}_2)^k$-action on a smooth connected closed manifold $M^n$
fixing a finite set can be associated to a $({\Bbb Z}_2)^k$-colored
regular graph. In this paper, we consider abstract graphs
$(\Gamma,\alpha)$ of $({\Bbb Z}_2)^k$-actions, called abstract
1-skeletons. We study when an abstract 1-skeleton is a colored graph
of some $({\Bbb Z}_2)^k$-action. We also study the existence of
faces of an abstract 1-skeleton (note that faces often have certain
geometric meanings if an abstract 1-skeleton is a colored graph of
some $({\Bbb Z}_2)^k$-action).
\end{abstract}

\section{Introduction}
Goresky, Kottwitz and MacPherson in \cite{gkm} showed that there is
an essential  connection between torus actions and regular graphs,
i.e., there is a large class of $T$-manifolds (called GKM
manifolds), each $M$ of which can be associated to a unique regular
graph $\Gamma_M$ so that the equivariant cohomology of $M$ can be
computed by the associated regular graph $\Gamma_M$, where $T$ is a
torus. A series of works by Guillemin and Zara further showed that
more topological and geometrical properties of $M$ can be read out
from $\Gamma_M$ (see, \cite{gz1}-\cite{gz4}).

\vskip .2cm

It was shown in \cite{l} that  the above idea can be extended to
2-torus actions, leading to study the equivariant cobordism
classification and the Smith problem for 2-torus actions.
Specifically, assume that $(\Phi, M)$ is an effective smooth
$G$-action on a smooth connected closed manifold $M$ fixing a finite
set (note that here $(\Phi, M)$ has less restriction than a GKM
manifold), where $G=({\Bbb Z}_2)^k$, a 2-torus of rank $k$. Then we
know from \cite[Section 2]{l} that the action $(\Phi, M)$ defines a
regular graph $\Gamma_{(\Phi, M)}$ with the vertex set $M^G$. This
graph is equipped with  a natural map (or a $G$-coloring) $\alpha$
from the set $E_{\Gamma_{(\Phi, M)}}$ of all edges of
$\Gamma_{(\Phi, M)}$ to all non-trivial elements of $\Hom(G,{\Bbb
Z}_2)$, satisfying  the following properties:
\begin{enumerate}
\item[(P1)] For each vertex $p$ of $\Gamma_{(\Phi, M)}$, the image
set $\alpha(E_p)$ spans $\Hom(G,{\Bbb Z}_2)$, where $E_p$ denotes
the set of all edges adjacent to $p$.

\item[(P2)] For each edge $e$ of $\Gamma_{(\Phi, M)}$,
$$\alpha(E_p)\equiv \alpha(E_q) \mod \alpha(e)$$
where $p$ and $q$ are two endpoints of $e$.
\end{enumerate} The pair $(\Gamma_{(\Phi,
M)}, \alpha)$ is called a {\em $G$-colored graph} of $(\Phi, M)$.
$(\Gamma_{(\Phi, M)}, \alpha)$ doesn't only contain the most
essential equivariant cobordism information of $(\Phi, M)$, but it
also indicates the relationship among $G$-representations on tangent
spaces at fixed points, so that it can be used to study the
equivariant cobordism classification and the Smith problem of
$(\Phi, M)$. Actually, from the colored graph $(\Gamma_{(\Phi,
  M)},\alpha)$ we  can read out all $G$-representations on
  tangent spaces at fixed points, which exactly consist of all images $\alpha(E_p), p\in
  M^G$. In particular,  $\{\alpha(E_p)| p\in
  M^G\}$  determines  a complete equivariant
  cobordism invariant  $\mathcal{P}_{(\Phi,M)}$ of
  $(\Phi, M)$, where  $\mathcal{P}_{(\Phi,M)}$ is obtained by
  deleting same pairs in $\{\alpha(E_p)| p\in
  M^G\}$.
 \vskip .2cm

 In this paper, we shall be concerned with an abstract pair
 $(\Gamma, \alpha)$, where $\Gamma$ is a finite $n$-valent regular graph without loops
 (i.e., edges with only an endpoint), and $\alpha$ satisfies (P1)
 and (P2). We call $(\Gamma, \alpha)$ {\em an abstract 1-skeleton of type ($n,k)$} (cf \cite[Definition 2.1.1]{gz2}).
 An easy observation shows that generally $(\Gamma, \alpha)$ may
 fail to be a colored graph induced by some $G$-action $(\Phi, M)$ (see also Example~\ref{ex1} of Section 2). A
 natural question is
\begin{enumerate}
\item[(Q1)] {\em When is $(\Gamma, \alpha)$ a  $G$-colored graph
of some $G$-action $(\Phi, M)$?}
\end{enumerate}
When $k=1$, we give a complete answer for (Q1), and when $k\geq 2$,
we give a partial answer for (Q1). We show that when $k\ge 2$, if
for each vertex $p\in V_\Gamma$, all vectors of $\alpha(E_p)$ are
pair-wise independent in $\Hom(G,{\Bbb Z}_2)$, then $\{\alpha(E_p)|
p\in
  V_\Gamma\}$ must be the fixed point data of some $G$-action $(\Phi,
  M)$. An example shows that the restriction of the pair-wise independence of $\alpha$ is necessary in the general case.
  However, we cannot make sure that $\Gamma$ is just a
  graph $\Gamma_{(\Phi,M)}$ of $(\Phi,
  M)$.

  \vskip .2cm

  We also consider those connected regular subgraphs of $(\Gamma, \alpha)$, each of
  which, with the restriction of $\alpha$ to it, is still an
  abstract 1-skeleton on its own right. We shall call those regular subgraphs with restrictions of $\alpha$
  the {\em faces} of
  $(\Gamma, \alpha)$. Generally, each face   has its geometric meaning  if $(\Gamma,
  \alpha)$ is a $G$-colored graph of a $G$-action $(\Phi, M)$.
  Actually it is often
  a $H$-colored graph of $H$-action on $N$, where $H$ is a
  subtorus of $G$, and $N$ is a component of the fixed point set
  of $G/H$ acting on $M$. We shall consider the existence of faces
  of an abstract 1-skeleton. Specifically, the following problem
  will be studied:
\begin{enumerate}
\item[(Q2)] {\em Let $(\Gamma, \alpha)$ be an abstract graph of
type $(k, n)$. Assume given a vertex $p\in V_\Gamma$ and a set of $m
(\geq 2)$ edges in $E_p$. Is there always a unique face containing
these $m$ edges?}
\end{enumerate}
The $l$-independence of $\alpha$ for $(Q2)$ is essentially
important. We shall show that if $\alpha$ is $l$-independent, then
 any $m$ edges in $E_p$  with $m<l$
  extend to a unique face. We shall also study the intersection
 property of faces, and then use it to consider the
 $n$-connectedness of $\Gamma$, obtaining a sufficient condition that
 $\Gamma$ is $n$-connected.

 \vskip .2cm
In the extreme case, an abstract 1-skeleton of type $(n,n)$ has its
own special properties. For example, the answer of (Q2) in this case
is always {\em yes}, so that  each such abstract graph
$(\Gamma,\alpha)$ can be
 associated with a unique simplicial poset, and so
 $(\Gamma,\alpha)$ has a geometric realization $|(\Gamma,\alpha)|$.
 If $(\Gamma, \alpha)$ is a colored graph of a small cover over a
 simple convex polytope $P$, then $|(\Gamma,\alpha)|$ is exactly the
 boundary of $P$. However, generally $|(\Gamma,\alpha)|$ even is not
 a closed manifold. In \cite{bl}, Bao and L\"u introduced the method of the skeletal
 expansion and gave a detailed investigation on $|(\Gamma,\alpha)|$.

\vskip .2cm

The paper is organized as follows. In Section 2, we introduce the
notions of an abstract 1-skeleton and the $l$-independence of
$\alpha$, and study the question (Q1). The notion of a face of
$(\Gamma, \alpha)$ is given in Section 3, and then the question (Q2)
is discussed. In Section 4, we consider the abstract 1-skeletons of
type $(n,n)$.

\vskip .2cm

The author would like to express his gratitude to Professor M.
Masuda for his valuable suggestions; especially for his help in the
proof of Proposition~\ref{n-con}. The author also would like to
express his gratitude to the referee, who did a careful reading. The
many suggestions and comments made by him or her considerably
improve the presentation of this paper.

\section{Abstract 1-skeletons}

Throughout the following one assumes that $G=({\Bbb Z}_2)^k$ with
$k\geq 1$.

\subsection{$G$-representations}

  Following \cite[Section 31]{cf}, let $R_n(G)$  denote the set generated by
the representation classes of dimension $n$, which naturally forms a
vector space over ${\Bbb Z}_2$. Then $R_*(G)=\sum_{n\geq 0}R_n(G)$
is a graded commutative algebra over ${\Bbb Z}_2$ with unit. The
multiplication in $R_*(G)$ is given by $[V_1]\cdot[V_2]=[V_1\oplus
V_2]$. Let $\Hom(G, {\Bbb Z}_2)$ be the set of all homomorphisms
$\rho: G\longrightarrow {\Bbb Z}_2=\{\pm 1\}$, which consists of
$2^k$ distinct homomorphisms, and let $\rho_0$ denote the trivial
element in $\Hom(G,{\Bbb Z}_2)$, i.e., $\rho_0(g)=1$ for all $g\in
G$. Every irreducible real representation of $G$ is one-dimensional
and has the form $\lambda_{\rho}: G\times {\Bbb R}\longrightarrow
{\Bbb R}$ with $\lambda_{\rho}(t, r) = \rho(t)\cdot r$ for some
$\rho\in \mbox{Hom}(G,{\Bbb Z}_2)$. Obviously there is a 1-1
correspondence between all irreducible real representations of $G$
and all elements of $\Hom(G,{\Bbb Z}_2)$. $\Hom(G,{\Bbb Z}_2)$ forms
an abelian group with addition given by
$(\rho+\sigma)(g)=\rho(g)\cdot\sigma(g)$, so it is also a vector
space over ${\Bbb Z}_2$ with standard basis $\{\rho_1,...,\rho_k\}$
where $\rho_i$ is defined by mapping $g=(g_1,...,g_i,...,g_k)$ to
$g_i$. Thus, we can identify $R_*(G)$ with the graded polynomial
algebra over ${\Bbb Z}_2$ generated by $\mbox{Hom}(G,{\Bbb Z}_2)$.
Namely,  $R_*(G)$ is isomorphic to the graded polynomial algebra
${\Bbb Z}_2[\rho_1,...,\rho_k]$.

\subsection{Abstract 1-skeletons} Let $\Gamma$ be a finite regular graph of
valence $n$ without loops such that $n\geq k$. If there is a map
$\alpha: E_\Gamma\longrightarrow \Hom(G,{\Bbb
Z}_2)\backslash{\{\rho_0\}}$ such that
\begin{enumerate}
\item[(1)] for each vertex $p\in V_\Gamma$, the image
$\alpha(E_p)$ spans $\Hom(G,{\Bbb Z}_2)$, and \item[(2)] for each
edge $e=pq\in E_\Gamma$,  $$\alpha(E_p)\equiv \alpha(E_q) \mod
\alpha(e)$$\end{enumerate} then the pair $(\Gamma, \alpha)$ is
called {\em an abstract 1-skeleton of type $(k,n)$}.

\vskip .2cm For each vertex $p\in V_\Gamma$, if any $l$ elements in
$\alpha(E_p)$ are linearly independent in $\Hom(G,{\Bbb Z}_2)$, then
one says that $\alpha$ is {\em $l$-independent} (cf \cite[Definition
2.1.2]{gz2}).

\begin{rem}
Since $\Hom(G,{\Bbb Z}_2)$ contains $2^k-1$ different nonzero
elements, one has that if  $\alpha$  is pairwise linearly
independent (i.e., 2-independent), then $k\geq 2$ and the valence
of $\Gamma$ is at most $2^k-1$.
\end{rem}

\subsection{Realization problem} Now let $(\Gamma, \alpha)$ be
 an abstract 1-skeleton of type $(k,n)$. We shall consider the question of whether  $(\Gamma,
 \alpha)$ can be realized as a $G$-colored graph of some $G$-action $(\Phi, M^n)$.

 \vskip .2cm

 When $k=1$, $\Hom(G,{\Bbb Z}_2)$ has only a non-trivial element,
 so for each vertex $p\in V_\Gamma$, all elements of $\alpha(E_p)$
 are same, and for any two $p, q\in V_\Gamma$,
 $\alpha(E_p)=\alpha(E_q)$.

 \begin{prop} \label{r1}
Let $(\Gamma, \alpha)$ is
 an abstract 1-skeleton of type $(1,n)$. Then $(\Gamma, \alpha)$
 is a $G$-colored graph of a $G$-action if and only if the
 number of vertices of $\Gamma$ is even.
 \end{prop}

 \begin{proof}
If $(\Gamma, \alpha)$
 is a $G$-colored graph of a $G$-action (i.e., an involution)
 $(\Phi, M)$, then the number of fixed points of $(\Phi, M)$ is
 the same as that of vertices of $\Gamma$. We know from \cite[Theorem 25.1]{cf}
 that if an involution fixes a finite set, then the number of
 fixed points must be even. Thus, the
 number of vertices of $\Gamma$ is even.

 \vskip .2cm

  Conversely, suppose that the  number of vertices of $\Gamma$ is
  even. With no loss of generality, assume that $\Gamma$ is
  connected. Given a $G$-action $\Phi$ on a connected closed manifold $M$ fixing a finite fixed set, we know from \cite[Section 2]{l} that
  any $n$-valent connected regular graph with the vertex set $M^G$ can be used as the colored graph of $(\Phi, M)$
  since all tangent $G$-representations at fixed points are the same. So, in order to complete the proof,
  it suffices to show that for each positive integer $l$, there is always an involution fixing $2l$ isolated points.
  In fact, a sphere $S^n$ always admits an involution fixing only two
  isolated points. Then the equivariant connected sum of  $l$
  copies of such involution on $S^n$ along their free orbits
  produces a new ${\Bbb Z}_2$-action $(\Psi, N)$, fixing $2l$ isolated points with same tangent $G$-representation.
   Thus, if the
 number of vertices of $\Gamma$ is even, then $(\Gamma, \alpha)$
 is a $G$-colored graph of a $G$-action.
 \end{proof}

 When $k\geq 2$, the problem becomes more complicated, but the 2-independence of $\alpha$ makes sure that
 $\{\alpha(E_p)| p\in V_\Gamma\}$ is the fixed point data of some
 $G$-action.

 \begin{prop}\label{r2}
Let $(\Gamma, \alpha)$ is
 an abstract 1-skeleton of type $(k,n)$ with $k\geq 2$. If $\alpha$ is
 2-independent, then $\{\alpha(E_p)| p\in V_\Gamma\}$ is the fixed point data of some
 $G$-action.
 \end{prop}

 \begin{proof}
According to tom Dieck-Kosniowski-Stong localization theorem (see
\cite{d}, \cite{ks} or \cite[Theorem 3.2]{l}), it suffices to show
that for any symmetric polynomial function $f(x_1,...,x_n)$ over
${\Bbb Z}_2$,
$$\sum_{p\in V_\Gamma}{{f(\alpha(E_p))}\over{\prod_{e\in
E_p}\alpha(e)}}\in {\Bbb Z}_2[\rho_1,...,\rho_k]$$ where
$f(\alpha(E_p))$ means that $x_1,...,x_n$ in $f(x_1,...,x_n)$ are
replaced by all elements in $\alpha(E_p)$. The idea of the following
proof is essentially due to Guillemin and Zara \cite[Theorem
2.2]{gz1}, but the proof is included here for local completeness.

\vskip .2cm If $\alpha$ is
 2-independent, then for each vertex $p$, all elements of
 $\alpha(E_p)$ are distinct. Thus, taking the common
denominator, $\sum_{p\in V_\Gamma}{{f(\alpha(E_p))}\over{\prod_{e\in
E_p}\alpha(e)}}$ becomes
$${h\over{\beta_1\beta_2\cdots \beta_u}}$$
where  $\beta_1,\beta_2,...,\beta_u\in \Hom(G,{\Bbb Z}_2)$ are
distinct. Now we want to show that $h$ is actually divisible by each
$\beta_i$. With no loss of generality, it suffices to prove that $h$
can be divided by $\beta_1$.  Let $V_1=\{p\in V_\Gamma| \beta_1\in
\alpha(E_p)\}$. It is easy to see that $V_1$ contains an even number
of vertices of $\Gamma$ since $\alpha$ is 2-independent. Take a
vertex $p$ in $V_1$; then there exists a unique edge $e$ in $E_p$
such that $\alpha(e)=\beta_1$. Let $q$ be another endpoint of $e$.
Since $\alpha(E_p)\equiv \alpha(E_q) \mod \alpha(e)=\beta_1$ and
$f(x_1,...,x_n)$ is symmetric, one has that $f(\alpha(E_p))\equiv
f(\alpha(E_q))\mod \beta_1$. Furthermore, by taking the common
denominator,  ${{f(\alpha(E_p))}\over {\prod_{y\in
E_p}\alpha(y)}}+{{f(\alpha(E_q))}\over {\prod_{z\in E_q}\alpha(z)}}$
becomes
$${{f(\alpha(E_p))g_1+f(\alpha(E_q))g_2}\over{\beta_1\beta'_2\cdots\beta'_v}}$$
where $g_1, g_2\in{\Bbb Z}_2[\rho_1,...,\rho_k]$.  Obviously,
$f(\alpha(E_p))g_1+f(\alpha(E_q))g_2$ can be divided by $\beta_1$.
Thus, we can write \begin{equation} \label{e1}
\sum_{p\in
V_1}{{f(\alpha(E_p))}\over{\prod_{e\in
E_p}\alpha(e)}}={{h_1}\over{\beta''_1\cdots \beta''_{u_1}}}
\end{equation} such that each $\beta''_i\in\Hom(G,{\Bbb Z}_2)$ is not
equal to $\beta_1$. Since $\beta_1\not\in\alpha(E_p)$ for $p\in
V_\Gamma\backslash V_1$, one has \begin{equation}
\label{e2}\sum_{p\in V_\Gamma\backslash
V_1}{{f(\alpha(E_p))}\over{\prod_{e\in
E_p}\alpha(e)}}={{h_2}\over{\beta'''_1\cdots \beta'''_{u_2}}}
\end{equation} such that each $\beta'''_i\in\Hom(G,{\Bbb Z}_2)$ is not
equal to $\beta_1$. Combining (\ref{e1}) and (\ref{e2}), one has
$$\sum_{p\in V_\Gamma}{{f(\alpha(E_p))}\over{\prod_{e\in
E_p}\alpha(e)}}={{h'}\over{\beta_2\cdots \beta_u}}.$$ This means
that $h$ is divisible by $\beta_1$.
 \end{proof}

The following example shows that generally, the restriction of
2-independence of $\alpha$ in Proposition~\ref{r2} is necessary.

\begin{exmp}\label{ex1}
{\em Figure~\ref{fig1} provides an abstract 1-skeleton of type $(3,
6)$ with four vertices $p,q,r,s$, which is not 2-independent.
\begin{figure}[h] \label{fig1}
    \input{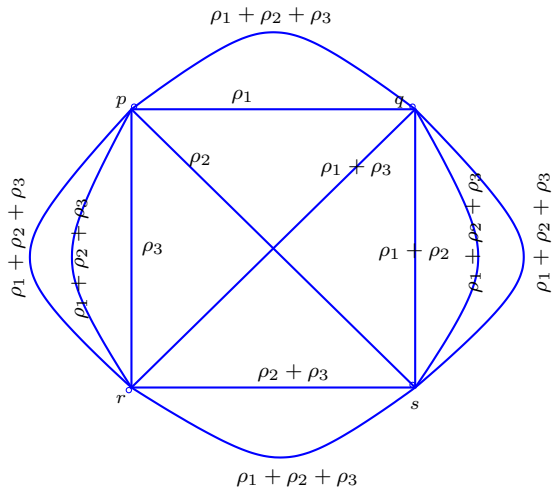}\centering
    \caption[a]{An abstract 1-skeleton of type $(3,6)$ with $\alpha$ no 2-independent.}
\end{figure}
 However, this abstract 1-skeleton is not a
$G$-colored graph of some $G$-action. This is because if one takes
$f(x_1,...,x_6)=\sigma_2(x_1,...,x_6)\sigma_3(x_1,...,x_6)$ with
$\deg f=5$, then by direct computations (cf \cite[Claim 2 and Remark
9]{l}),
$$\sum_{u\in\{p,q,r,s\}}{{f(\alpha(E_u))}\over{\prod_{e\in E_u}\alpha(e)}}\not=0$$
so $\sum_{u\in\{p,q,r,s\}}{{f(\alpha(E_u))}\over{\prod_{e\in
E_u}\alpha(e)}}\not\in{\Bbb Z}_2[\rho_1,\rho_2,\rho_3]$.  Note that
for each $u$, $\deg \prod_{e\in E_u}\alpha(e)=6$.  }
\end{exmp}

\begin{rem}
We tried to show that  if $\alpha$ is
 2-independent, then $(\Gamma,\alpha)$ is  a $G$-colored graph of some
 $G$-action, but failed. Even so, it is extremely tempting to
 {\em conjecture} that this is true.
\end{rem}
\section{Faces}

 Suppose that $(\Gamma, \alpha)$ is an
abstract 1-skeleton of type $(k,n)$.  Let $\Gamma'$ be a connected
$\ell$-valent subgraph of $\Gamma$ where $0\leq \ell\leq n$.

\vskip .2cm

 We say that $(\Gamma',
\alpha|_{\Gamma'})$ is a {\em $\ell$-dimensional face} of
$(\Gamma, \alpha)$ if there is a subspace $K$ of $\Hom(G,{\Bbb
Z}_2)$ such that
\begin{enumerate}
\item[(1)] for each vertex $p\in V_{\Gamma'}$, the image
$\alpha(E_p|_{\Gamma'})$ spans $K$, and \item[(2)] for each edge
$e=pq\in E_{\Gamma'}$,  $$\alpha(E_p|_{\Gamma'})\equiv
\alpha(E_q|_{\Gamma'}) \mod \alpha(e).$$\end{enumerate}

\vskip .2cm

 Obviously,  each edge $e$ of $(\Gamma,\alpha)$, with the restriction of $\alpha$ to it, is a 1-dimensional
  face. $(\Gamma, \alpha)$ is
the union of some $n$-faces, and in particular, $(\Gamma,\alpha)$
itself is a unique $n$-face if $\Gamma$ is connected. Note that each
vertex of $\Gamma$ is a $0$-face.  The following example shows that
generally, for $m$ edges of $E_p$ at a vertex $p$ where $1<m<n$,
there may  be no $m$-face containing the $m$ edges.

\vskip .2cm

\begin{exmp} \label{face}
 {\em Figure~\ref{fig2} gives an abstract 1-skeleton $(\Gamma, \alpha)$ of type $(3, 4)$ such that $\alpha$ is
  three-independent, but there are three edges at the vertex $p$
  colored by $\rho_1, \rho_2, \rho_3$, respectively
such that they cannot extend to a 3-face.}
\begin{figure}[h] \label{fig2}
    \input{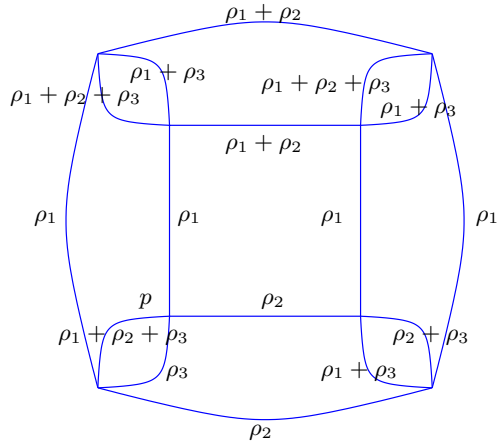}\centering
    \caption[a]{An abstract 1-skeleton of type $(3,4)$ with $\alpha$  three-independent.}
\end{figure}
\end{exmp}
\vskip .2cm

The existence of a face of $(\Gamma,\alpha)$ depends heavily on
the $l$-independence of $\alpha$. Our result for (Q2) is stated as
follows.

\begin{prop}\label{f1}
Suppose that $(\Gamma, \alpha)$ is an abstract 1-skeleton of type
$(k, n)$ such that  $\alpha$  is $l$-independent where $1<l\leq k$.
Then any $m (<l)$ edges $e_1, ..., e_m$ of $E_p$ at $p\in V_\Gamma$
can extend to a unique $m$-face. In particular, when $k=n$, any $m
(\leq n)$ edges of $E_p$ at $p\in V_\Gamma$ can always extend to a
unique $m$-face.
\end{prop}

\begin{proof}
It is trivial when $l=2$. If $l\geq 3$, write
$E=\{e_1,...,e_m\}\subset E_p$, one then proceeds as follows. One
begins with $p$ as being the starting point. Since $\alpha$ is
$l$-independent and $m<l$, the space $K=\Span\{\alpha(e)|e\in E\}$
has dimension $m$ and is different from the space
$\Span\{\alpha(e)|e\in E'\}$, where $E'\subset E_p$ contains $m$
edges and $E'\not= E$.  Take an edge $\bar{e}$ in $E$ and let $q$ be
another endpoint of $\bar{e}$, one then has that there is a unique
subset $F$ of containing $m$ edges in $E_q$ such that
$K=\Span\{\alpha(e)|e\in F\}$ and $\alpha(E)\equiv \alpha(F)\mod
\alpha(\bar{e})$. One may further carry out this procedure at $q$ as
follows: take an edge $\bar{e}'$ in $F$ and let $r$ be another
vertex of $\bar{e}'$, then one finds a unique subset $H$ of $E_r$ of
containing $m$ edges such that $K=\Span\{\alpha(e)|e\in H\}$ and
$\alpha(F)\equiv \alpha(H)\mod \alpha(\bar{e}')$.
 Continuing this procedure, since $\Gamma$ is assumed
to be finite, finally one can obtain a unique connected $m$-valent
subgraph containing $e_1,...,e_m$ as desired.
\end{proof}

\begin{rem}
As shown in \cite[Proposition 2.1.3]{gz2}, the three-independence of
$\alpha$ for a GKM graph $(\Gamma, \alpha)$ determines a unique
connection $\theta$. This is also ture for our abstract graphs in
mod 2 category. However, the three-independence of $\alpha$ cannot
determine the existence of faces of dimension more than 2,  as is
shown in Example~\ref{face}.
\end{rem}

We also can further obverse the  property of the intersection of
faces of $(\Gamma, \alpha)$.

\begin{prop}\label{f2}
Suppose that $(\Gamma, \alpha)$ is an abstract 1-skeleton of type
$(k,n)$ such that  $\alpha$  is $l$-independent, where $l\leq k$.
Then for faces $F^{m_1},...,F^{m_s}$ of dimension less than $l$ in
$(\Gamma, \alpha)$ with $s>1$, their intersection
$F^{m_1}\cap\cdots\cap F^{m_s}$ is either empty or the disjoint
union of lower-dimensional faces.
\end{prop}
\begin{proof}
Suppose that the intersection $F^{m_1}\cap\cdots\cap
F^{m_s}\not=\emptyset$. Then $F^{m_1}\cap\cdots\cap F^{m_s}$
contains at least one vertex (i.e., a $0$-face) $p$ of $\Gamma$.
Suppose $E_p\vert_{F^{m_1}}\cap\cdots\cap E_p\vert_{F^{m_s}}$
contains only $m$ edges in $E_p$. Then  $0\leq m\leq \min(m_1,..,
m_s)$ and  by Proposition~\ref{f1} the intersection
$F^{m_1}\cap\cdots\cap F^{m_s}$ must contain a $m$-face $F^m$
containing the $m$ edges in $E_p$. If there is also another vertex
$q\in F^{m_1}\cap\cdots\cap F^{m_s}$ such that $q\not\in F^m$, in
the above way one may obtain a face $F^{m'}$ containing $q$.

\vskip .2cm One claims that  $F^{m'}$ is another connected component
different from $F^m$ in the intersection $F^{m_1}\cap\cdots\cap
F^{m_s}$. With no loss of generality, one assumes that $m'\leq m$.
If the dimension of $F^{m'}$ is zero (i.e., $F^{m'}=\{q\}$), then
obviously the claim holds. If the dimension of $F^{m'}$ is more than
zero and if the intersection of $F^{m'}$ and $F^m$ is nonempty, then
there is a vertex $o$ in $\Gamma$ such that $o\in F^{m'}\cap F^m$.
Obviously, all edges of $F^{m'}$ with   endpoint $o$ cannot belong
to $F^m$. Indeed, otherwise $F^{m'}$ would be a subface of $F^m$, so
$q\in F^m$, which gives a contradiction. Therefore, there is an edge
$e$
 in $F^{m'}$ with   endpoint $o$ such that $e\not\in F^m$. Since
$F^m\subset F^{m_1}\cap\cdots\cap F^{m_s}$ and $F^{m'}\subset
F^{m_1}\cap\cdots\cap F^{m_s}$, one has that $ F^{m_1}\cap\cdots\cap
F^{m_s}$ contains at least $m+1$ edges at $E_o$ such that these
$m+1$ edges extend a $(m+1)$-face $F^{m+1}$ with $F^m$ as its a
subface, and in particular, $F^{m+1}$ contains the vertex $p$. Since
$F^{m+1}\subset F^{m_1}\cap\cdots\cap F^{m_s}$, this means that
$F^{m_1}\cap\cdots\cap F^{m_s}$ contains $m+1$ edges of $E_p$. This
is a contradiction.
\end{proof}

The following is an application for the $d$-connectedness of a
graph. Before one discusses this, let us review the notion for the
$d$-connectedness of a graph and the Whitney Theorem (cf \cite{g}).

\vskip .2cm Let $\Gamma$ be a graph. A {\em path} $\sigma$ with
two endpoints $a, b$ in $\Gamma$ is a subgraph of $\Gamma$ having
as vertices the vertices $v_0=a, v_1,...,v_{l-1}, v_l=b$ of
$\Gamma$ and as edges the edges $v_{i-1}v_i, i=1,...,l$ of
$\Gamma$. Two paths $\sigma_1$ and $\sigma_2$  with common
endpoints $a, b$ are called {\em disjoint} provided the
intersection $\sigma_1\cap \sigma_2$ consists of $a$ and $b$ only.
A graph $\Gamma$ is {\em connected} provided for each pair of its
vertices there is a path in $\Gamma$ having these vertices as
endpoints. A graph $\Gamma$ is {\em $d$-connected} provided for
every pair of vertices of $\Gamma$ there exist $d$ pairwise
disjoint paths in $\Gamma$ having these vertices as endpoints.

\vskip .3cm

\begin{thm} {\em (Whitney)}  A graph $\Gamma$ with at least $d+1$ vertices
is $d$-connected if and only if every subgraph of $\Gamma$ obtained
by deleting from $\Gamma$ any $d-1$ or fewer vertices and the edges
incident to them is connected.
\end{thm}

\begin{lem} \label{2-con}
Suppose that $(\Gamma, \alpha)$ is an abstract 1-skeleton of type $(k,n)$
with $\Gamma$ connected such that  $\alpha$ is at least three
independent. Then $\Gamma$ must be at least 2-connected.
\end{lem}
\begin{proof}
If $\Gamma$ is 1-connected but not 2-connected, then by Whitney
Theorem there must be at least an edge $e=pq$ in $\Gamma$ such that
$\Gamma\backslash\{e\}$ is disconnected. Take an edge $e'$ in
$E_p\backslash\{e\}$, since $\alpha$ is at least three independent,
by Proposition~\ref{f1} there should be a 2-face $F$ containing $e$
and $e'$. However, this is impossible since $\Gamma\backslash\{e\}$
is disconnected.
\end{proof}

\begin{prop}\label{n-con}
Suppose that $(\Gamma, \alpha)$ is an abstract 1-skeleton of type
$(k,n)$ with $\Gamma$ connected such that $\alpha$ is at least three
independent. If the intersection of any two faces of dimension less
than $3$ in $(\Gamma, \alpha)$ is either connected or empty, then
$\Gamma$ is $n$-connected.
\end{prop}

\begin{proof}
Suppose that $\Gamma$ is not $n$-connected. Then by Whitney
Theorem and Lemma~\ref{2-con}, there is at least one subset $S$
with its cardinality $2\leq \vert S\vert\leq n-1$ of the vertex
set $V_\Gamma$ such that $\Gamma\backslash S$ is disconnected,
where $\Gamma\backslash S$ denotes the subgraph of $\Gamma$
obtained by removing  all vertices in $S$ and all edges adjacent
to those vertices. Let $S_{\min}$ be {\em minimal} (i.e., its
cardinality is minimal) among subsets $S$ for which
$\Gamma\backslash S$ is disconnected. Fix a vertex $p$ of
$S_{\min}$.

\vskip .2cm

{\em Claim 1.} Each connected component of $\Gamma\backslash
S_{\min}$ contains a vertex adjacent to~$p$.

\vskip .2cm

 Take any vertex $q$ of $\Gamma\backslash S_{\min}$. Since
$S_{\min}$ is minimal, $\Gamma\backslash (S_{\min}\backslash\{p\})$
is connected; so there is a path in $\Gamma\backslash
(S_{\min}\backslash\{p\})$ from $q$ to $p$. Cutting the last edge
from this path produces a path in $\Gamma\backslash S_{\min}$ from
$q$ to a vertex adjacent to $p$. This proves the claim 1.

\vskip .2cm Let $q, q'$ be vertices of $\Gamma\backslash S_{\min}$
which belong to different connected components of
$\Gamma\backslash S_{\min}$ and are adjacent to $p$. Then one has
that

\vskip .2cm

{\em Claim 2.} The 2-face determined by the two edges joining $p$
and $q, q'$ contains a vertex in $S_{\min}$ except $p$.

\vskip .2cm

 Indeed, if the 2-face does not contain any vertex in $S_{\min}$ except
$p$, then it gives a path joining $q$ and $q'$ in
$\Gamma\backslash S_{\min}$, but this contradicts the assumption
that $q$ and $q'$ belong to different connected component of
$\Gamma\backslash S_{\min}$.

\vskip .2cm

Now let $r$ be the number of vertices in $S_{\min}$ adjacent to $p$.
Then one has

\vskip .2cm

{\em Claim 3.} There are at least $n-r-1$  two-faces which contain
$p$ and some other vertex of $S_{\min}$ not adjacent to $p$.

\vskip .2cm

 Let $N(p)$ be the set of vertices of $\Gamma\backslash S_{\min}$
adjacent to $p$. The cardinality of $N(p)$ is $n-r$. A connected
component of $N(p)$ means the intersection of $N(p)$ with a
connected component of $\Gamma\backslash S_{\min}$. By Claim 1, the
number of connected components of $N(p)$ is the same as that of
connected components of $\Gamma\backslash S_{\min}$. Suppose that a
connected component $C$ of $N(p)$ has $d$ vertices, where $1\leq
d\leq n-r-1$. Choosing one vertex from $C$ and the other one from
$N(p)\backslash C$, one gets a 2-face $F$ containing two edges
joining $p$ and the chosen two vertices. $F$ contains $p$ and a
vertex in $S_{\min}$ except $p$, say $v$, by Claim 2. If $v$ is
adjacent to $p$, then the intersection of $F$ and a 1-face joining
$p$ and $v$ consists of $p, v$ (i.e., two 0-faces) and hence is
disconnected. But this impossible, so $v$ is not adjacent to $p$.
There are $d(n-r-d)$  such 2-faces and it is easy to see that
$d(n-r-d)\geq n-r-1$.

\vskip .2cm

On the other hand,   the number of vertices in
$S_{\min}\backslash\{p\}$ which are not adjacent to $p$ is $\vert
S_{\min}\vert-r-1\leq n-r-2$. This together with Claim 3 implies
that there are at least two different 2-faces $F_1, F_2$ which
contain $p$ and some other vertex, say $q$, of $S_{\min}$ not
adjacent to $p$. Thus, the intersection of two 2-faces $F_1, F_2$
contains $p$ and $q$, which are not adjacent to each other. Since
$\alpha$ is at least three independent, by Proposition~\ref{f2} the
intersection $F_1\cap F_2$ consists of the disjoint union of some
1-faces and 0-faces. This means that the intersection $F_1\cap F_2$
is disconnected. This is a contradiction. Therefore $\vert
S_{\min}\vert\geq n$ and $\Gamma$ is $n$-connected.
\end{proof}

\vskip .3cm

\begin{figure}[h]\label{fig3}
    \input{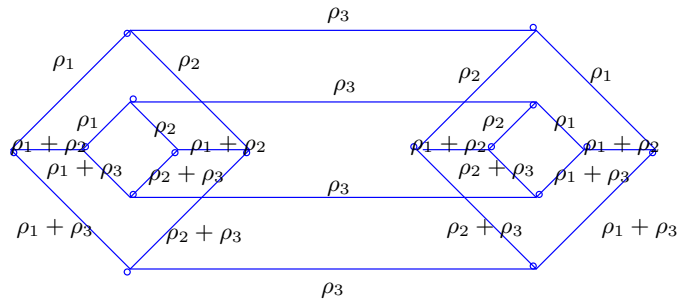}\centering
    \caption[a]{ An abstract 1-skeleton of type $(3,3)$}\label{}
\end{figure}

{\em Note.} The inverse of Proposition~\ref{n-con} is generally
untrue. For example, Figure~\ref{fig3}  gives an abstract 1-skeleton
$(\Gamma, \alpha)$ of type $(3,3)$ such that $\alpha$ is
three-independent. Obviously, $\Gamma$ is 3-connected, but a direct
observation shows that there are at least two 2-faces in
$(\Gamma,\alpha)$ such that their intersection is disconnected.

\section{Abstract 1-skeletons of type $(n,n)$}

In this section, we are concerned with  the case $k=n$. Suppose that
$(\Gamma,\alpha)$ is an  abstract 1-skeleton  of type $(n,n)$. In
this case, $\alpha$ is $n$-independent so by Proposition~\ref{f1}
any $m (\leq n)$ edges of $E_p$ at $p\in V_\Gamma$ can extend to a
unique $m$-face. Let $\mathcal{F}_{(\Gamma, \alpha)}$ denote the set
of all faces of an abstract 1-skeleton $(\Gamma,\alpha)$ of type
$(n,n)$.  We call an $(n-1)$-face of $\mathcal{F}_{(\Gamma,
\alpha)}$ a {\em facet} (cf \cite{mmp}).

\begin{lem} \label{p1}
Let $F^m$ be a $m$-face with $m<n$ in $\mathcal{F}_{(\Gamma,
\alpha)}$. Then there  exist $n-m$ facets $F_1,...,F_{n-m}$ in
$\mathcal{F}_{(\Gamma, \alpha)}$ such that $F^m$ is a connected
component of the intersection $F_1\cap\cdots\cap F_{n-m}$.
\end{lem}
\begin{proof}
For a vertex $p$ in $F^m$, $F^m$ contains $m$ edges of $E_p$. Write
$E_p=\{e_1,...,e_n\}$, with no loss one assumes that $e_1,...,e_m$
belong to $F^m$. By Proposition~\ref{f1},  any $n-1$ edges of $E_p$
determine a facet, so  there are exactly $n$ facets with $p$ as a
vertex. Obviously there exactly are $n-m$ facets $F_1,...,F_{n-m}$
with $p$ as a vertex such that each $F_i$ contains these $m$ edges
$e_1,...,e_m$, and  does not contain at least  one edge  from
$E_p\backslash\{e_1,...,e_m\}$.  Since any facet of containing
$e_1,...,e_m$ must contain $F^m$, one has that  each facet $F_i$
contains $F^m$. The lemma then follows from Proposition~\ref{f2}.
\end{proof}

Every face of $(\Gamma, \alpha)$ corresponds to  a subspace of
$\Hom(G,{\Bbb Z}_2)$.
\begin{lem}\label{p2}
Suppose that  $(\Gamma, \alpha)$ is an abstract 1-skeleton of type
$(n,n)$. Let $\Gamma'$ be a connected $m$-valent subgraph of
$\Gamma$. Then $(\Gamma', \alpha\vert_{E_{\Gamma'}})$ is an $m$-face
with $m>0$ if and only if $\alpha(E_{\Gamma'})$ spans an
$m$-dimensional  of $\Hom(G, {\Bbb Z}_2)$.
\end{lem}
\begin{proof}
If $(\Gamma', \alpha\vert_{E_{\Gamma'}})$ is an $m$-face, then for
any two different vertices $p_1, p_2$ in $V_{\Gamma'}$,
$\alpha(E_{\Gamma'}\vert_{p_1})$ and
$\alpha(E_{\Gamma'}\vert_{p_2})$ span the same $m$-dimensional
vector subspace of $\Hom(G, {\Bbb Z}_2)$, so $\alpha(E_{\Gamma'})$
spans an $m$-dimensional  of $\Hom(G, {\Bbb Z}_2)$.

\vskip .2cm

Conversely, suppose  $\alpha(E_{\Gamma'})$ spans a $m$-dimensional
of $\Hom(G, {\Bbb Z}_2)$. Then, since $\alpha$ is $n$-independent,
for any two vertices $p_1$ and $p_2$ in $V_{\Gamma'}$,
$\Span\alpha(E_{\Gamma'}\vert_{p_1})$ and
$\Span\alpha(E_{\Gamma'}\vert_{p_2})$ must be the same
$m$-dimensional subspace of $\Hom(G, {\Bbb Z}_2)$. An easy
observation also shows that for each edge $e=pq$ in $\Gamma'$,
$\alpha(E_{\Gamma'}\vert_{p})\equiv \alpha(E_{\Gamma'}\vert_{q})\mod
\alpha(e)$. Thus,  $(\Gamma', \alpha\vert_{E_{\Gamma'}})$ is an
$m$-face.
\end{proof}

As vector spaces over ${\Bbb Z}_2$, $\Hom(G,{\Bbb Z}_2)$ and
$\Hom({\Bbb Z}_2, G)$ are dual to each other. By Lemma~\ref{p2},
each $m$-face $F^m=(\Gamma', \alpha\vert_{E_{\Gamma'}})$ with $m>0$
of $\mathcal{F}_{(\Gamma, \alpha)}$ actually corresponds to a unique
$(n-m)$-dimensional subspace $J$ of $\Hom({\Bbb Z}_2, G)$ such that
for each $a^*\in J$ and each $a\in \alpha(E_{\Gamma'})$,
$$a^*(a)=0.$$
  Thus, each facet  of
$\mathcal{F}_{(\Gamma, \alpha)}$ corresponds to a unique nonzero
element of $\Hom({\Bbb Z}_2, G)$. This gives a map
$$\lambda:\overline{\mathcal{F}}_{(\Gamma,
\alpha)}\longrightarrow \Hom({\Bbb Z}_2, G)$$ where
$\overline{\mathcal{F}}_{(\Gamma, \alpha)}$ denotes the set of all
facets of $\mathcal{F}_{(\Gamma, \alpha)}$. For each vertex $p$, let
$F_1,...,F_n$ be $n$ facets with $p$ as a vertex. Then it is easy to
see that $\lambda(F_1),...,\lambda(F_n)$ are linearly independent in
$\Hom({\Bbb Z}_2, G)$. Given a $m$-face $F^m$, by Lemma~\ref{p1} one
has that there  exist $n-m$ facets $F_1,...,F_{n-m}$ in
$\mathcal{F}_{(\Gamma, \alpha)}$ such that $F^m$ is a connected
component of the intersection $F_1\cap\cdots\cap F_{n-m}$.
 Then $F^m$
corresponds to the $(n-m)$-dimensional vector space spanned by
$\lambda(F_1),...,\lambda(F_{n-m})$. Here one calls  $\lambda$ the
{\em characteristic function}. Combining the above arguments, one
has

\begin{prop}\label{dual}
Suppose that  $(\Gamma, \alpha)$ is an abstract 1-skeleton of type
$(n,n)$. Then  $\alpha$ and $\lambda$  determine each other.
\end{prop}

\begin{rem}
If $(\Gamma, \alpha)$ is a colored graph of a small cover $M$ over a
simple convex polytope $P$, then $\Gamma$ is exactly the 1-skeleton
of $P$, and each facet of $\mathcal{F}_{(\Gamma, \alpha)}$
corresponds to a facet of $P$. For the notion of a small cover, see
\cite{dj}. Thus, the map $\lambda$ defined above is actually the
characteristic function of the small cover $M$.  Furthermore,
combining the GKM theory and the Davis-Januszkiewicz theory for
small covers together, we see that the face ring of $P$ over ${\Bbb
Z}_2$ is isomorphic to
$$\{f: V_\Gamma\longrightarrow {\Bbb Z}_2[\rho_1,...,\rho_n]\big|
f(p)\equiv f(q) \mod \alpha(e) \text{ for } e\in E_p\cap E_q\}.$$ It
should be pointed out that an abstract 1-skeleton of type $(n,n)$ is
an analogue of a torus graph introduced by Maeda, Masuda and Panov
in \cite{mmp}.  They have shown that  the equivariant cohomology of
a torus graph is isomorphic to the face ring of the associated
simplicial poset. This is a generalization of the above isomorphism.
\end{rem}

Now suppose that  $(\Gamma, \alpha)$ is an abstract 1-skeleton of
type $(n,n)$ with $\Gamma$  connected.  Then $\mathcal{F}_{(\Gamma,
\alpha)}$ contains a unique $n$-face $(\Gamma, \alpha)$. As an
analogue,
 $(\Gamma, \alpha)$  has many
same properties as a torus graph, which are stated as follows:
\begin{enumerate}
\item[(i)]$\mathcal{F}_{(\Gamma, \alpha)}$ forms a simplicial poset
of rank $n$ with respect to  reversed inclusion with $(\Gamma,
\alpha)$ as smallest element, denoted by $\mathcal{P}_{(\Gamma,
\alpha)}$.
\item[(ii)]$\mathcal{P}_{(\Gamma, \alpha)}$ is a face poset of a simplicial
complex $K$ if and only if all possible non-empty intersections of
facets of $\mathcal{F}_{(\Gamma, \alpha)}$ are connected (see also
\cite[Proposition 5.1]{mmp}).
\end{enumerate}

As a consequence of Lemma~\ref{p1},  Proposition~\ref{n-con} and
Property (ii), one has that

\begin{cor}
If $\mathcal{P}_{(\Gamma, \alpha)}$ is a face poset of a simplicial
complex $K$, then $\Gamma$ is $n$-connected.
\end{cor}

We know from \cite{s} or \cite{mmp} that as a simplicial poset,
$\mathcal{P}_{(\Gamma, \alpha)}$ determines a regular CW-complex
${\Bbb K}_{\mathcal{P}_{(\Gamma, \alpha)}}$, which is
$(n-1)$-dimensional. By $\vert(\Gamma,\alpha)\vert$ one denotes the
underlying space of this cell complex, and one calls it {\em the
geometric realization of $(\Gamma,\alpha)$}.

\begin{rem} The
 geometric realization $\vert(\Gamma, \alpha)\vert$ of $(\Gamma, \alpha)$ has a direct connection with the topology of manifolds.
  As an interesting
topic, the study of $\vert(\Gamma, \alpha)\vert$ has been carried on
independently in \cite{bl}. For example, it can be shown that each
closed combinatorial manifold can be realizable by $\vert(\Gamma,
\alpha)\vert$,
 and that the geometric realization
 of an abstract 1-skeleton $(\Gamma, \alpha)$ of type $(4,4)$ is a closed 3-manifold if and only if the $f$-vector of
 $\mathcal{P}_{(\Gamma,\alpha)}$ satisfies $f_1=f_0+f_3$, where $f_i$
denote the number of $a\in \mathcal{P}_{(\Gamma,\alpha)}$ for which
the segment $[\hat{0},a]$ is a boolean algebra of rank $i+1$.
\end{rem}


\begin{thebibliography}{99}

\bibitem[BL]{bl} Z. Q. Bao and Z. L\"u, {\em Manifolds associated with $({\Bbb Z}_2)^n$-colored regular
graphs},  arXiv: math/0609557.
\bibitem[CF]{cf} P.E. Conner and E.E. Floyd, {\em Differentiable
 periodic maps}, Ergebnisse Math. Grenzgebiete, N. F., Bd. {\bf 33},
 Springer-Verlag, Berlin, 1964.
 \bibitem[D]{d} T. tom Dieck, {\em Characteristic numbers of
 $G$-manifolds. I}, Invent. Math. {\bf 13} (1971), 213-224.
 \bibitem[DJ]{dj} M. Davis and T. Januszkiewicz, {\em Convex
polytope, Coxeter orbifolds and torus actions}, Duke Math. J. {\bf
62} (1991), 417-451.
\bibitem[G]{g} B. Gr\"unbaum,  {\em Convex Polytopes}, Graduate Texts in Mathematics, {\bf 221}, Springer,
2003.
\bibitem[GKM]{gkm} M. Goresky, R. Kottwitz, and R. MacPherson, {\em Equivariant cohomology, Koszul duality, and
   the localization theorem},  Invent. Math.  {\bf 131} (1998), 25-83.
   \bibitem[GZ1]{gz1} V. Guillemin and C. Zara, {\em
  Equivariant de Rham theory and graphs}, Asian J. Math. {\bf 3}
  (1999), 47-76.
\bibitem[GZ2]{gz2} V. Guillemin and C. Zara, {\em 1-Skeleta, Betti numbers, and equivariant cohomology},
 Duke Math. J. {\bf 107} (2001), 283-349.
  \bibitem[GZ3]{gz3} V. Guillemin and C. Zara, {\em The existence of generating families for the cohomology ring of
  a graph},  Advances in Mathematics {\bf 174} (2003), 115-153.
  \bibitem[GZ4]{gz4} V. Guillemin and C. Zara, {\em $G$-actions on graphs},
    Internat. Math. Res. Notices  {\bf 10} (2001), 519--542.
     \bibitem[KS]{ks} C. Kosniowski and R.E. Stong, {\em $({\Bbb Z}_2)^k$-actions
 and characteristic numbers}, Indiana Univ. Math. J. {\bf 28} (1979),
 723-743.
    \bibitem[L]{l} Z. L\"u, {\em Graphs and $({\Bbb Z}_2)^k$-actions},  arXiv:math/0508643.
    \bibitem[MMP]{mmp} H. Maeda, M. Masuda and T. Panov, {\em Torus
    graphs and simplicial posets}, Adv. Math. {\bf 212} (2007), 458-483.
        \bibitem[S]{s} R.P. Stanley, {\em $f$-vectors and $h$-vectors of
 simplicial posets}, Journal of Pure and Applied Algebra {\bf 71}
 (1991), 319-331.



 \end{thebibliography}
\end{document}